\documentclass[11pt]{article}
\usepackage{epsfig}

\def\be{\begin{equation}}
\def\ee{\end{equation}}
\def\bea{\begin{eqnarray}}
\def\eea{\end{eqnarray}}
\def\bes{\begin{eqnarray*}}
\def\ees{\end{eqnarray*}}

\def\nn{\nonumber}
\def\<{\langle}
\def\>{\rangle}
\def\lb{\label}
\def\bs{\setminus}

\def\R{{\bf R}}
\def\C{{\bf C}}
\def\Z{{\bf Z}}
\def\N{{\bf N}}
\def\U{{\bf U}}
\def\Q{{\bf Q}}
\def\T{{\bf T}}

\def\ga{{\gamma}}

\def\th{{\theta}}
\def\om{{\omega}}
\def\Om{{\Omega}}
\def\ep{{\epsilon}}
\def\lm{{\lambda}}

\def\sg{{\sigma}}

\def\Sg{{\Sigma}}
\def\vf{{\varphi}}

\def\H{{\cal H}}
\def\T{{\cal T}}

\def\P{{\cal P}}

\def\Nn{{\cal N}}

\def\Sp{{\rm Sp}}

\def\dm{{\rm \diamond}}

\def\ol#1{\overline{#1}}  

\def\hb{\vrule height0.18cm width0.14cm $\,$}

\def\ol#1{\overline{#1}}  

\def\mapright#1{\smash{\mathop{\longrightarrow}\limits^{#1}}}

\title{Existence of closed characteristics on compact \\convex  hypersurfaces in $\R^{2n}$}
\author{Wei Wang\thanks{Partially supported by National Natural Science Foundation of China No.10801002, China Postdoctoral Science Foundation No.200801021, Foundation for the Author of National Excellent Doctoral Dissertation of PR China No. 201017.
E-mail: alexanderweiwang@yahoo.com.cn, wangwei@math.pku.edu.cn  }\\
Beijing International Center for  Mathematical  Research\\
Key Laboratory of Pure and Applied Mathematics\\
School of Mathematical Science \\ Peking University, Beijing 100871 \\
PEOPLES REPUBLIC OF CHINA \\ }
\date{ Dec.  16th, 2011}

\begin{document}

\maketitle

\begin{abstract}
{\it  In this paper, we prove there exist at least $[\frac{n+1}{2}]+1$ geometrically distinct 
closed characteristics on every compact convex hypersurface $\Sg$ in $\R^{2n}$. Moreover,
there exist at least $[\frac{n}{2}]+1$ geometrically distinct non-hyperbolic closed characteristics on  $\Sg$ in $\R^{2n}$ provided the number of  geometrically distinct closed characteristics on $\Sg$ is finite.}

\end{abstract}



{\bf Running title}: Closed characteristics on convex hypersurfaces

\renewcommand{\theequation}{\thesection.\arabic{equation}}
\renewcommand{\thefigure}{\thesection.\arabic{figure}}

\setcounter{equation}{0}
\section{Introduction and main results}

Let $\Sigma$ be a fixed $C^3$ compact convex hypersurface in
$\R^{2n}$, i.e., $\Sigma$ is the boundary of a compact and strictly
convex region $U$ in $\R^{2n}$. We denote the set of all such
hypersurfaces by $\H(2n)$. Without loss of generality, we suppose
$U$ contains the origin. We consider closed characteristics $(\tau,
y)$ on $\Sigma$, which are solutions of the following problem \be
\left\{\matrix{\dot{y}=JN_\Sigma(y), \cr
               y(\tau)=y(0), \cr }\right. \lb{1.1}\ee
where $J=\left(\matrix{0 &-I_n\cr
                I_n  & 0\cr}\right)$,
$I_n$ is the identity matrix in $\R^n$, $\tau>0$, $N_\Sigma(y)$ is
the outward normal vector of $\Sigma$ at $y$ normalized by the
condition $N_\Sigma(y)\cdot y=1$. Here $a\cdot b$ denotes the
standard inner product of $a, b\in\R^{2n}$.
A closed characteristic
$(\tau, y)$ is {\it prime}, if $\tau$ is the minimal period of $y$.
Two closed characteristics $(\tau, y)$ and $(\sigma, z)$ are {\it
geometrically distinct},  if $y(\R)\not= z(\R)$. We denote by
$\T(\Sg)$ the set of all geometrically distinct closed
characteristics on $\Sg$. A closed characteristic $(\tau,y)$ is {\it
non-degenerate}, if $1$ is a Floquet multiplier of $y$ of precisely
algebraic multiplicity $2$, and is {\it elliptic}, if all the
Floquet multipliers of $y$ are on ${\bf U}=\{z\in\C\,|\,|z|=1\}$,
i.e., the unit circle in the complex plane. It is  {\it hyperbolic},
if $1$ is a double Floquet multiplier of it and all the other
Floquet multipliers of $y$ are away from  ${\bf
U}=\{z\in\C\,|\,|z|=1\}$.

It is surprising enough that A.M. Liapunov in \cite{Lia1} of 1892
and J. Horn in \cite{Hor1} of 1903 were able to prove the following
great result: {\it Suppose $H: \R^{2n}\rightarrow\R$ is analytic,
$\sigma(JH^{\prime\prime}(0))=\{\pm\sqrt{-1}\omega_1, \ldots
,\pm\sqrt{-1}\omega_n\}$ are purly imaginary and satisfy
$\frac{\omega_i}{\omega_j}\notin\Z$ for all $i, j$. Then there
exists $\epsilon_0>0$ so small that }
 \be
^{\#}\T(H^{-1}(\epsilon))\ge n, \qquad \forall \; 0< \epsilon\le
\epsilon_0.\lb{1.2}\ee 
This deep result was greatly improved by A.
Weinstein in \cite{Wei1} of 1973. He was able to prove that for
$H\in C^2(\R^{2n}, \R)$, if $H^{\prime\prime}(0)$ is positive
definite, then there exists $\epsilon_0>0$ small such that
(\ref{1.2}) still holds. In \cite{EL}, I. Ekeland and  J. Lasry
proved that if there exists $x_0\in\R^{2n}$ such that
\bea r\le |x-x_0|\le R,\qquad \forall x\in\Sg
\nn\eea
and $\frac{R}{r}<\sqrt{2}$.
Then $\Sg$ carries at least $n$ geometrically distinct
closed characteristics.

Note that we have the following example of weakly non-resonant
ellipsoid: Let $r=(r_1,\ldots, r_n)$ with $r_i>0$ for $1\le i\le n$.
Define \bea \mathcal{E}_n(r)=\left\{z=(x_1, \ldots,x_n,
y_1,\ldots,y_n)\in\R^{2n}\left |\frac{}{}\right.
\frac{1}{2}\sum_{i=1}^n\frac{x_i^2+y_i^2}{r_i^2}=1\right\}\nn\eea¡¡
where $\frac{r_i}{r_j}\notin\Q$ whenever $i\neq j$. In this case,
the corresponding Hamiltonian system is linear and all the solutions
can be computed explicitly. Thus it is easy to verify that
$^{\#}\T(\mathcal{E}_n(r)) = n$ and all the closed characteristics
on $\mathcal{E}_n(r)$ are elliptic and non-degenerate, i.e. its
linearized Poincar\'e map splits into $n-1$ two dimensional rotation
matrix $\left(\matrix{\cos\theta &\sin\theta\cr -\sin\theta & \cos\theta\cr}\right)$
with $\frac{\theta}{\pi}\notin\Q$ and one $\left(\matrix{1 &1\cr 0 & 1\cr}\right)$.

Based on the above facts, there is a long standing conjecture on the
number of closed characteristics on compact convex hypersurfaces in
$\R^{2n}$: \be \,^{\#}\T(\Sg)\ge n, \qquad \forall \; \Sg\in\H(2n).
\lb{1.3}\ee
Since the pioneering works \cite{Rab1} of P. Rabinowitz and
\cite{Wei2} of A. Weinstein in 1978 on the existence of at least one
closed characteristic on every hypersurface in $\H(2n)$, the
existence of multiple closed characteristics on $\Sg\in\H(2n)$ has
been deeply studied by many mathematicians.
When $n\ge 2$,  in 1987-1988 I. Ekeland-L.
Lassoued, I. Ekeland-H. Hofer, and A, Szulkin (cf. \cite{EkL1},
\cite{EkH1}, \cite{Szu1}) proved
$$ \,^{\#}\T(\Sg)\ge 2, \qquad \forall\,\Sg\in\H(2n). $$
 In \cite{HWZ1} of 1998, H. Hofer-K. Wysocki-E. Zehnder
proved that $\,^{\#}\T(\Sg)=2$ or $\infty$ holds for every
$\Sg\in\H(4)$. In \cite{LoZ1} of 2002, Y. Long and C. Zhu further proved
$$ \;^{\#}\T(\Sg)\ge [\frac{n}{2}]+1, \qquad \forall\, \Sg\in \H(2n), $$
where we denote by $[a]\equiv\max\{k\in\Z\,|\,k\le a\}$. 
 In \cite{WHL} of 2007,
W. Wang, X. Hu and Y. Long proved $\,^{\#}\T(\Sg)\ge 3$ for
every $\Sg\in\H(6)$, which gave a confirmed answer to the above
conjecture for $n=3$.

There are some related results considering the stability problem, in
\cite{Eke2} of Ekeland in 1986 and \cite{Lon2} of Long in 1998, for
any $\Sg\in\H(2n)$ the existence of at least one non-hyperbolic
closed characteristic on $\Sg$ was proved provided
$^\#\T(\Sg)<+\infty$. Ekeland proved also in \cite{Eke2} the
existence of at least one elliptic closed characteristic on $\Sg$
provided $\Sg\in\H(2n)$ is $\sqrt{2}$-pinched. In \cite{DDE1} of
1992, Dell'Antonio, D'Onofrio and Ekeland proved the existence of at
least one elliptic closed characteristic on $\Sg$ provided
$\Sg\in\H(2n)$ satisfies $\Sg=-\Sg$.  
 In \cite{Lon3} of 2000, Y. Long proved  that
$\Sg\in\H(4)$ and $\,^{\#}\T(\Sg)=2$ imply that both of the
closed characteristics must be elliptic. In \cite{WHL}, the authors
proved further that $\Sg\in\H(4)$ and $\,^{\#}\T(\Sg)=2$ imply
that both of the closed characteristics must be irrationally
elliptic.In \cite{LoZ1} of 2002, Long
and Zhu  proved when $^\#\T(\Sg)<+\infty$, there exists
at least one elliptic closed characteristic and there are at least
$[\frac{n}{2}]$ geometrically distinct closed characteristics on
$\Sg$ possessing irrational mean indices. which are then
non-hyperbolic. In \cite{W1}, the author proved that if
$^\#\T(\Sg)<+\infty$, then at least two closed characteristics
on $\Sg$ must possess irrational mean indices; if
$^\#\T(\Sg)=3$, then there are at least two elliptic closed
characteristics on $\Sg$, where $\Sg\in\H(6)$.

Motivated by these results, we prove the following results in this paper.

{\bf Theorem 1.1.} {\it   There exist at least $[\frac{n+1}{2}]+1$ geometrically distinct 
closed characteristics on every compact convex hypersurface $\Sg$  in $\R^{2n}$, i.e., we have
$^{\#}\T(\Sg)\ge [\frac{n+1}{2}]+1$ for any $\Sg\in\H(2n)$.}

{\bf Theorem 1.2.} {\it There exist at least $[\frac{n}{2}]+1$ geometrically distinct non-hyperbolic closed characteristics on  $\Sg\in\H(2n)$ provided the number of  geometrically distinct closed characteristics on $\Sg$ is finite.}

The proof of the main Theorems is given in Section 4. Mainly ingredients
in the proof include: the critical point theory for closed
characteristics established in \cite{WHL}, Morse theory and
the index iteration theory developed by Long and his coworkers,
specially the common index jump theorem of Long and Zhu (Theorem 4.3
of \cite{LoZ1}, cf. Theorem 11.2.1 of \cite{Lon4}). In Section 2, we
review briefly the equivariant Morse theory for closed
characteristics on compact convex hypersurfaces in $\R^{2n}$
developed in \cite{WHL}. In Section 3, we review the 
index iteration theory developed by Long and his coworkers.

In this paper, let $\N$, $\N_0$, $\Z$, $\Q$, $\R$, and $\R^+$ denote
the sets of natural integers, non-negative integers, integers,
rational numbers, real numbers, and positive real numbers
respectively. Denote by $a\cdot b$ and $|a|$ the standard inner
product and norm in $\R^{2n}$. Denote by $\langle\cdot,\cdot\rangle$
and $\|\cdot\|$ the standard $L^2$-inner product and $L^2$-norm. For
an $S^1$-space $X$, we denote by $X_{S^1}$ the homotopy quotient of
$X$ module the $S^1$-action, i.e., $X_{S^1}=S^\infty\times_{S^1}X$.
We define the functions \be \left\{\matrix{[a]=\max\{k\in\Z\,|\,k\le
a\}, & E(a)=\min\{k\in\Z\,|\,k\ge a\} , \cr
                   \varphi(a)=E(a)-[a],   \cr}\right. \lb{1.4}\ee
Specially, $\varphi(a)=0$ if $ a\in\Z\,$, and $\varphi(a)=1$ if
$a\notin\Z\,$. In this paper we use only $\Q$-coefficients for all
homological modules. For a $\Z_m$-space pair $(A, B)$, let
$H_{\ast}(A, B)^{\pm\Z_m}= \{\sigma\in H_{\ast}(A,
B)\,|\,L_{\ast}\sigma=\pm \sigma\}$, where $L$ is a generator of the
$\Z_m$-action.

\setcounter{equation}{0}
\section{ Critical point theory for closed characteristics}

In the rest of this paper, we fix a $\Sg\in\H(2n)$ and assume the
following condition on $\Sg$:

\noindent (F) {\bf There exist only finitely many geometrically
distinct closed characteristics \\$\quad \{(\tau_j, y_j)\}_{1\le
j\le k}$ on $\Sigma$. }

In this section, we review briefly the equivariant Morse theory for
closed characteristics on $\Sg$ developed in \cite{WHL} and
\cite{W1} which will be needed in Section 4 of this paper. All the
details of proofs can be found in \cite{WHL} and \cite{W1}.

Let $\hat{\tau}=\inf\{\tau_j|\;1\le j\le k\}$. Note that here
$\tau_j$'s are prime periods of $y_j$'s for $1\le j\le k$. Then by
\S2 of \cite{WHL}, for any $a>\hat{\tau}$, we can construct a
function $\varphi_a\in C^\infty(\R,\R^+)$ which has $0$ as its
unique critical point in $[0,\,+\infty)$ such that $\varphi_a$ is
strictly convex for $t\ge 0$. Moreover,
$\frac{\varphi_a^\prime(t)}{t}$ is strictly decreasing for $t> 0$
together with $\lim_{t\rightarrow
0^+}\frac{\varphi_a^\prime(t)}{t}=1$ and
$\varphi_a(0)=0=\varphi_a^\prime(0)$. More precisely, we define
$\varphi_a$ via Propositions 2.2 and 2.4 in \cite{WHL}. The precise
dependence of $\varphi_a$ on $a$ is explained in Remark 2.3 of
\cite{WHL}.

Define the Hamiltonian function $H_a(x)=a\varphi_a(j(x))$ and
consider the fixed period problem \be
\left\{\matrix{\dot{x}(t)=JH_a^\prime(x(t)), \cr
     x(1)=x(0).         \cr }\right. \lb{2.1}\ee
Then  $H_a\in C^3(\R^{2n}\setminus\{0\}, \R)\cap C^1(\R^{2n}, \R)$
is strictly convex. Solutions of (\ref{2.1}) are $x\equiv0$ and
$x=\rho y(\tau t)$ with
$\frac{\varphi_a^\prime(\rho)}{\rho}=\frac{\tau}{a}$, where $(\tau,
y)$ is a solution of (\ref{1.1}). In particular, nonzero solutions
of (\ref{2.1}) are one to one correspondent to solutions of
(\ref{1.1}) with period $\tau<a$.

In the following, we use the Clarke-Ekeland dual action principle.
As usual, let $G_a$ be the Fenchel transform of $H_a$ defined by
$G_a(y)=\sup\{x\cdot y-H_a(x)\;|\; x\in \R^{2n}\}$. Then $G_a\in
C^2(\R^{2n}\bs\{0\},\R)\cap C^1(\R^{2n},\R)$ is strictly convex. Let
\be L_0^2(S^1, \;\R^{2n})= \left\{u\in L^2([0, 1],\;\R^{2n})
   \left|\frac{}{}\right.\int_0^1u(t)dt=0\right\}.  \lb{2.2}\ee
Define a linear operator $M: L_0^2(S^1,\R^{2n})\to
L_0^2(S^1,\R^{2n})$ by $\frac{d}{dt}Mu(t)=u(t)$,
$\int_0^1Mu(t)dt=0$. The dual action functional on $L_0^2(S^1,
\;\R^{2n})$ is defined by \be
\Psi_a(u)=\int_0^1\left(\frac{1}{2}Ju\cdot Mu+G_a(-Ju)\right)dt.
   \lb{2.3}\ee
Then the functional $\Psi_a\in C^{1, 1}(L_0^2(S^1,\; \R^{2n}),\;\R)$
is bounded from below and satisfies the Palais-Smale condition.
Suppose $x$ is a solution of (\ref{2.1}). Then $u=\dot{x}$ is a
critical point of $\Psi_a$. Conversely, suppose $u$ is a critical
point of $\Psi_a$. Then there exists a unique $\xi\in\R^{2n}$ such
that $Mu-\xi$ is a solution of (\ref{2.1}). In particular, solutions
of (\ref{2.1}) are in one to one correspondence with critical points
of $\Psi_a$. Moreover, $\Psi_a(u)<0$ for every critical point
$u\not= 0$ of $\Psi_a$.

Suppose $u$ is a nonzero critical point of $\Psi_a$. Then following
\cite{Eke3} the formal Hessian of $\Psi_a$ at $u$ is defined by
$$ Q_a(v,\; v)=\int_0^1 (Jv\cdot Mv+G_a^{\prime\prime}(-Ju)Jv\cdot Jv)dt, $$
which defines an orthogonal splitting $L_0^2=E_-\oplus E_0\oplus
E_+$ of $L_0^2(S^1,\; \R^{2n})$ into negative, zero and positive
subspaces. The index of $u$ is defined by $i(u)=\dim E_-$ and the
nullity of $u$ is defined by $\nu(u)=\dim E_0$. Let $u=\dot{x}$ be
the critical point of $\Psi_a$ such that $x$ corresponds to the
closed characteristic $(\tau,\,y)$ on $\Sigma$. Then the index
$i(u)$ and the nullity $\nu(u)$ defined above coincide with the
Ekeland indices defined by I. Ekeland in \cite{Eke1} and
\cite{Eke3}. Specially $1\le \nu(u)\le 2n-1$ always holds.

We have a natural $S^1$-action on $L_0^2(S^1,\; \R^{2n})$ defined by
$\th\cdot u(t)=u(\th+t)$ for all $\th\in S^1$ and $t\in\R$. Clearly
$\Psi_a$ is $S^1$-invariant. For any $\kappa\in\R$, we denote by \be
\Lambda_a^\kappa=\{u\in L_0^2(S^1,\;
\R^{2n})\;|\;\Psi_a(u)\le\kappa\}.
          \lb{2.4}\ee
For a critical point $u$ of $\Psi_a$, we denote by \be
\Lambda_a(u)=\Lambda_a^{\Psi_a(u)}
  =\{w\in L_0^2(S^1,\; \R^{2n}) \;|\; \Psi_a(w)\le\Psi_a(u)\}.\lb{2.5}\ee
Clearly, both sets are $S^1$-invariant. Since the $S^1$-action
preserves $\Psi_a$, if $u$ is a critical point of $\Psi_a$, then the
whole orbit $S^1\cdot u$ is formed by critical points of $\Psi_a$.
Denote by $crit(\Psi_a)$ the set of critical points of $\Psi_a$.
Note that by the condition (F), the number of critical orbits of
$\Psi_a$ is finite. Hence as usual we can make the following
definition.

{\bf Definition 2.1.} {\it Suppose $u$ is a nonzero critical point
of $\Psi_a$ and $\Nn$ is an $S^1$-invariant open neighborhood of
$S^1\cdot u$ such that $crit(\Psi_a)\cap(\Lambda_a(u)\cap
\Nn)=S^1\cdot u$. Then the $S^1$-critical modules of $S^1\cdot u$
are defined by}
$$ C_{S^1,\; q}(\Psi_a, \;S^1\cdot u)
=H_{q}((\Lambda_a(u)\cap\Nn)_{S^1},\; ((\Lambda_a(u)\setminus
S^1\cdot u)\cap\Nn)_{S^1}). $$

We have the following proposition for critical modules.

{\bf Proposition 2.2.} (Proposition 3.2 of \cite{WHL}) {\it The
critical module $C_{S^1,\;q}(\Psi_a, \;S^1\cdot u)$ is independent
of $a$ in the sense that if $x_i$ are solutions of (\ref{2.1}) with
Hamiltonian functions $H_{a_i}(x)\equiv a_i\varphi_{a_i}(j(x))$ for
$i=1$ and $2$ respectively such that both $x_1$ and $x_2$ correspond
to the same closed characteristic $(\tau, y)$ on $\Sigma$. Then we
have}
$$ C_{S^1,\; q}(\Psi_{a_1}, \;S^1\cdot\dot {x}_1) \cong
  C_{S^1,\; q}(\Psi_{a_2}, \;S^1\cdot \dot {x}_2), \quad \forall q\in \Z. $$

Now let $u\neq 0$ be a critical point of $\Psi_a$ with multiplicity
$mul(u)=m$, i.e., $u$ corresponds to a closed characteristic
$(m\tau, y)\subset\Sigma$ with $(\tau, y)$ being prime. Hence
$u(t+\frac{1}{m})=u(t)$ holds for all $t\in \R$ and the orbit of
$u$, namely, $S^1\cdot u\cong S^1/\Z_m\cong S^1$. Let $f: N(S^1\cdot
u)\rightarrow S^1\cdot u$ be the normal bundle of $S^1\cdot u$ in
$L_0^2(S^1,\; \R^{2n})$ and let $f^{-1}(\theta\cdot u)=N(\theta\cdot
u)$ be the fibre over $\theta\cdot u$, where $\theta\in S^1$. Let
$DN(S^1\cdot u)$ be the $\varrho$-disk bundle of $N(S^1\cdot u)$ for
some $\varrho>0$ sufficiently small, i.e., $DN(S^1\cdot u)=\{\xi\in
N(S^1\cdot u)\;| \; \|\xi\|<\varrho\}$ and let $DN(\theta\cdot
u)=f^{-1}(\th\cdot u)\cap DN(S^1\cdot u)$ be the disk over
$\theta\cdot u$. Clearly, $DN(\theta\cdot u)$ is $\Z_m$-invariant
and we have $DN(S^1\cdot u)=DN(u)\times_{\Z_m}S^1$, where the
$Z_m$-action is given by
$$ (\th, v, t)\in \Z_m\times DN(u)\times S^1\mapsto
        (\th\cdot v, \;\theta^{-1}t)\in DN(u)\times S^1. $$
Hence for an $S^1$-invariant subset $\Gamma$ of $DN(S^1\cdot u)$, we
have $\Gamma/S^1=(\Gamma_u\times_{\Z_m}S^1)/S^1=\Gamma_u/\Z_m$,
where $\Gamma_u=\Gamma\cap DN(u)$. Since $\Psi_a$ is not $C^2$ on
$L_0^2(S^1,\; \R^{2n})$, we need to use a finite dimensional
approximation introduced by Ekeland in order to apply Morse theory.
More precisely, we can construct a finite dimensional submanifold
$\Gamma(\iota)$ of $L_0^2(S^1,\; \R^{2n})$ which admits a
$\Z_\iota$-action with $m|\iota$. Moreover $\Psi_a$ and
$\Psi_a|_{\Gamma(\iota)}$ have the same critical points.
$\Psi_a|_{\Gamma(\iota)}$ is $C^2$ in a small tubular neighborhood
of the critical orbit $S^1\cdot u$ and the Morse index and nullity
of its critical points coincide with those of the corresponding
critical points of $\Psi_a$.  Let \be D_\iota N(S^1\cdot
u)=DN(S^1\cdot u)\cap\Gamma(\iota), \quad D_\iota N(\theta\cdot
u)=DN(\theta\cdot u)\cap\Gamma(\iota). \lb{2.6}\ee Then we have \be
C_{S^1,\; \ast}(\Psi_a, \;S^1\cdot u) \cong H_\ast(\Lambda_a(u)\cap
D_\iota N(u),\;
    (\Lambda_a(u)\setminus\{u\})\cap D_\iota N(u))^{\Z_m}. \lb{2.7}\ee
Now we can apply the results of Gromoll and Meyer in \cite{GrM1} to
the manifold $D_{p\iota}N(u^p)$ with $u^p$ as its unique critical
point, where $p\in\N$. Then $mul(u^p)=pm$ is the multiplicity of
$u^p$ and the isotropy group $\Z_{pm}\subseteq S^1$ of $u^p$ acts on
$D_{p\iota}N(u^p)$ by isometries. According to Lemma 1 of
\cite{GrM1}, we have a $\Z_{pm}$-invariant decomposition of
$T_{u^p}(D_{p\iota}N(u^p))$
$$ T_{u^p}(D_{p\iota}N(u^p))
=V^+\oplus V^-\oplus V^0=\{(x_+, x_-, x_0)\}  $$ with $\dim
V^-=i(u^p)$, $\dim V^0=\nu(u^p)-1$ and a $\Z_{pm}$-invariant
neighborhood $B=B_+\times B_-\times B_0$ for $0$ in
$T_{u^p}(D_{p\iota}N(u^p))$ together with two $Z_{pm}$-invariant
diffeomorphisms
$$\Phi :B=B_+\times B_-\times B_0\rightarrow
\Phi(B_+\times B_-\times B_0)\subset D_{p\iota}N(u^p)$$ and
$$ \eta : B_0\rightarrow W(u^p)\equiv\eta(B_0)\subset D_{p\iota}N(u^p)$$
such that $\Phi(0)=\eta(0)=u^p$ and \be
\Psi_a\circ\Phi(x_+,x_-,x_0)=|x_+|^2 - |x_-|^2 +
\Psi_a\circ\eta(x_0),
    \lb{2.8}\ee
with $d(\Psi_a\circ \eta)(0)=d^2(\Psi_a\circ\eta)(0)=0$. As
\cite{GrM1}, we call $W(u^p)$ a local {\it characteristic manifold}
and $U(u^p)=B_-$ a local {\it negative disk} at $u^p$. By the proof
of Lemma 1 of \cite{GrM1}, $W(u^p)$ and $U(u^p)$ are
$\Z_{pm}$-invariant. Then we have \bea && H_\ast(\Lambda_a(u^p)\cap
D_{p\iota}N(u^p),\;
  (\Lambda_a(u^p)\setminus\{u^p\})\cap D_{p\iota}N(u^p)) \nn\\
&&\qquad = H_\ast (U(u^p),\;U(u^p)\setminus\{u^p\}) \otimes
H_\ast(W(u^p)\cap \Lambda_a(u^p),\; (W(u^p)\setminus\{u^p\})\cap
\Lambda_a(u^p)),
  \lb{2.9}\eea
where \be H_q(U(u^p),U(u^p)\setminus\{u^p\} )
    = \left\{\matrix{\Q, & {\rm if\;}q=i(u^p),  \cr
                      0, & {\rm otherwise}. \cr}\right.  \lb{2.10}\ee
Now we have the following proposition.

{\bf Proposition 2.3.} (Proposition 3.10 of \cite{WHL}) {\it Let
$u\neq 0$ be a critical point of $\Psi_a$ with $mul(u)=1$. Then for
all $p\in\N$ and $q\in\Z$, we have \be C_{S^1,\; q}(\Psi_a,
\;S^1\cdot u^p)\cong \left(\frac{}{}H_{q-i(u^p)}(W(u^p)\cap
\Lambda_a(u^p),\; (W(u^p)\setminus\{u^p\})\cap
\Lambda_a(u^p))\right)^{\beta(u^p)\Z_p},
  \lb{2.11}\ee
where $\beta(u^p)=(-1)^{i(u^p)-i(u)}$. Thus \be C_{S^1,\; q}(\Psi_a,
\;S^1\cdot u^p)=0, \quad {\rm for}\;\;
   q<i(u^p) \;\;{\rm or}\;\;q>i(u^p)+\nu(u^p)-1. \lb{2.12}\ee
In particular, if $u^p$ is non-degenerate, i.e., $\nu(u^p)=1$, then}
\be C_{S^1,\; q}(\Psi_a, \;S^1\cdot u^p)
    = \left\{\matrix{\Q, & {\rm if\;}q=i(u^p)\;{\rm and\;}\beta(u^p)=1,  \cr
                      0, & {\rm otherwise}. \cr}\right.  \lb{2.13}\ee

We make the following definition.

{\bf Definition 2.4.} {\it Let $u\neq 0$ be a critical point of
$\Psi_a$ with $mul(u)=1$. Then for all $p\in\N$ and $l\in\Z$, let
\bea k_{l, \pm 1}(u^p)&=&\dim\left(\frac{}{}H_l(W(u^p)\cap
\Lambda_a(u^p),\;
(W(u^p)\setminus\{u^p\})\cap \Lambda_a(u^p))\right)^{\pm\Z_p}, \nn\\
k_l(u^p)&=&\dim\left(\frac{}{}H_l(W(u^p)\cap \Lambda_a(u^p),
(W(u^p)\setminus\{u^p\})\cap
\Lambda_a(u^p))\right)^{\beta(u^p)\Z_p}. \nn\eea $k_l(u^p)$'s are
called critical type numbers of $u^p$. }

We have the following properties for critical type numbers.

{\bf Proposition 2.5.} (Proposition 3.13 of \cite{WHL}) {\it Let
$u\neq 0$ be a critical point of $\Psi_a$ with $mul(u)=1$. Then
there exists a minimal $K(u)\in \N$ such that
$$ \nu(u^{p+K(u)})=\nu(u^p),\quad i(u^{p+K(u)})-i(u^p)\in 2\Z, $$
and $k_l(u^{p+K(u)})=k_l(u^p)$ for all $p\in \N$ and $l\in\Z$. We
call $K(u)$ the minimal period of critical modules of iterations of
the functional $\Psi_a$ at $u$.   }

For a closed characteristic $(\tau,y)$ on $\Sigma$, we denote by
$y^m\equiv (m\tau, y)$ the $m$-th iteration of $y$ for $m\in\N$. Let
$a>\tau$ and choose $\vf_a$ as above. Determine $\rho$ uniquely by
$\frac{\vf_a'(\rho)}{\rho}=\frac{\tau}{a}$. Let $x=\rho y(\tau t)$
and $u=\dot{x}$. Then we define the index $i(y^m)$ and nullity
$\nu(y^m)$ of $(m\tau,y)$ for $m\in\N$ by
$$ i(y^m)=i(u^m), \qquad \nu(y^m)=\nu(u^m). $$
These indices are independent of $a$ when $a$ tends to infinity. Now
the mean index of $(\tau,y)$ is defined by
$$ \hat{i}(y)=\lim_{m\rightarrow\infty}\frac{i(y^m)}{m}. $$
Note that $\hat{i}(y)>2$ always holds which was proved by Ekeland
and Hofer in \cite{EkH1} of 1987 (cf. Corollary 8.3.2 and Lemma
15.3.2 of \cite{Lon4} for a different proof).

By Proposition 2.2, we can define the critical type numbers
$k_l(y^m)$ of $y^m$ to be $k_l(u^m)$, where $u^m$ is the critical
point of $\Psi_a$ corresponding to $y^m$. We also define
$K(y)=K(u)$. Then we have the following.

{\bf Proposition 2.6.} {\it We have $k_l(y^m)=0$ for $l\notin [0,
\nu(y^m)-1]$ and it can take only values $0$ or $1$ when $l=0$ or
$l=\nu(y^m)-1$. Moreover, the following properties hold (cf. Lemma
3.10 of \cite{BaL1}, \cite{Cha1} and \cite{MaW1}):

(i) $k_0(y^m)=1$ implies $k_l(y^m)=0$ for $1\le l\le \nu(y^m)-1$.

(ii) $k_{\nu(y^m)-1}(y^m)=1$ implies $k_l(y^m)=0$ for $0\le l\le
\nu(y^m)-2$.

(iii) $k_l(y^m)\ge 1$ for some $1\le l\le \nu(y^m)-2$ implies
$k_0(y^m)=k_{\nu(y^m)-1}(y^m)=0$.

(iv) If $\nu(y^m)\le 3$, then at most one of the $k_l(y^m)$'s for
$0\le l\le \nu(y^m)-1$ can be non-zero.

(v) If $i(y^m)-i(y)\in 2\Z+1$ for some $m\in\N$, then $k_0(y^m)=0$.}

{\bf Proof.} By Definition 2.4 we have
$$ k_l(y^m)\le \dim H_l(W(u^m)\cap \Lambda_a(u^m),\;
(W(u^m)\setminus\{u^m\})\cap \Lambda_a(u^m))\equiv \eta_l(y^m). $$
Then from Corollary 1.5.1 of \cite{Cha1} or Corollary 8.4 of
\cite{MaW1}, (i)-(iv) hold.

For (v), if $\eta_0(y^m)=0$, then (v) follows directly from
Definition 2.4.

By Corollary 8.4 of \cite{MaW1}, $\eta_0(y^m)=1$ if and only if
$u^m$ is a local minimum in the local characteristic manifold
$W(u^m)$. Hence $(W(u^m)\cap
\Lambda_a(u^m),\;(W(u^m)\setminus\{u^m\})\cap
\Lambda_a(u^m))=(\{u^m\},\; \emptyset)$. By Definition 2.4, we have:
\bea k_{0, +1}(u^m) &=& \dim H_0(W(u^m)\cap \Lambda_a(u^m),\;
    (W(u^m)\setminus\{u^m\})\cap \Lambda_a(u^m))^{+\Z_m}\nn\\
&=& \dim H_0(\{u^m\})^{+\Z_m}\nn\\
&=& 1.   \nn\eea This implies $k_0(u^m)=k_{0, -1}(u^m)=0$. \hfill\hb

Let $\Psi_a$ be the functional defined by (\ref{2.3}) for some
$a\in\R$ large enough and let $\varepsilon>0$ be small enough such
that $[-\varepsilon, +\infty)\setminus\{0\}$ contains no critical
values of $\Psi_a$. Denote by $I_a$ the greatest integer in $\N_0$
such that $I_a<i(\tau, y)$ hold for all closed characteristics
$(\tau,\, y)$ on $\Sigma$ with $\tau\ge a$. Then by Section 5 of
\cite{WHL}, we have \be H_{S^1,\; q}(\Lambda_a^{-\varepsilon} )
\cong H_{S^1,\; q}( \Lambda_a^\infty)
  \cong H_q(CP^\infty), \quad \forall q<I_a.  \lb{2.14}\ee
For any $q\in\Z$, let \be  M_q(\Lambda_a^{-\varepsilon})
  =\sum_{1\le j\le k,\,1\le m_j<a/\tau_j} \dim C_{S^1,\;q}(\Psi_a, \;S^1\cdot u_j^{m_j}).
  \lb{2.15} \ee
Then the equivariant Morse inequalities for the space
$\Lambda_a^{-\varepsilon}$ yield \bea M_q(\Lambda_a^{-\varepsilon})
       &\ge& b_q(\Lambda_a^{-\varepsilon}),\lb{2.16}\\
M_q(\Lambda_a^{-\varepsilon}) &-& M_{q-1}(\Lambda_a^{-\varepsilon})
    + \cdots +(-1)^{q}M_0(\Lambda_a^{-\varepsilon}) \nn\\
&\ge& b_q(\Lambda_a^{-\varepsilon}) -
b_{q-1}(\Lambda_a^{-\varepsilon})
   + \cdots + (-1)^{q}b_0(\Lambda_a^{-\varepsilon}), \lb{2.17}\eea
where $b_q(\Lambda_a^{-\varepsilon})=\dim H_{S^1,\;
q}(\Lambda_a^{-\varepsilon})$. Now we have the following Morse
inequalities for closed characteristics.

{\bf Theorem 2.7.} {\it Let $\Sigma\in \H(2n)$ satisfy
$\,^{\#}\T(\Sg)<+\infty$. Denote all the geometrically
distinct closed characteristics by $\{(\tau_j,\; y_j)\}_{1\le j\le
k}$. Let \bea
M_q&=&\lim_{a\rightarrow+\infty}M_q(\Lambda_a^{-\varepsilon}),\quad
                  \forall q\in\Z,\lb{2.18}\\
b_q &=& \lim_{a\rightarrow+\infty}b_q(\Lambda_a^{-\varepsilon})=
\left\{\matrix{1, & {\rm if\;}q\in 2\N_0,  \cr
                      0, & {\rm otherwise}. \cr}\right.  \lb{2.19}
\eea Then we have}
\bea M_q &\ge& b_q,\lb{2.20}\\
 M_q-M_{q-1}+\cdots +(-1)^{q}M_0 &\ge& b_q-b_{q-1}+\cdots +(-1)^{q}b_0,
   \qquad\forall \;q\in\Z. \lb{2.21}\eea

{\bf Proof.} As we have mentioned before, $\hat i(y_j)>2$ holds for
$1\le j\le k$. Hence the Ekeland index satisfies
$i(y_j^m)=i(u_j^m)\to\infty$ as $m\to\infty$ for $1\le j\le k$. Note
that $I_a\to +\infty$ as $a\to +\infty$. Now fix a $q\in\Z$ and a
sufficiently great $a>0$. By Propositions 2.2, 2.3 and (\ref{2.15}),
$M_i(\Lambda_a^{-\varepsilon})$ is invariant for all $a>A_q$ and
$0\le i\le q$, where $A_q>0$ is some constant. Hence (\ref{2.18}) is
meaningful. Now for any $a$ such that $I_a>q$,
(\ref{2.14})-(\ref{2.17}) imply that (\ref{2.19})-(\ref{2.21}) hold.
\hfill\hb

Recall that for a principal $U(1)$-bundle $E\to B$, the Fadell-Rabinowitz index
(cf. \cite{FaR1}) of $E$ is defined to be $\sup\{k\;|\, c_1(E)^{k-1}\not= 0\}$,
where $c_1(E)\in H^2(B,\Q)$ is the first rational Chern class. For a $U(1)$-space,
i.e., a topological space $X$ with a $U(1)$-action, the Fadell-Rabinowitz index is
defined to be the index of the bundle $X\times S^{\infty}\to X\times_{U(1)}S^{\infty}$,
where $S^{\infty}\to CP^{\infty}$ is the universal $U(1)$-bundle.

As in P.199 of \cite{Eke3}, choose some $\alpha\in(1,\, 2)$ and associate with $U$
a convex function $H$ such that $H(\lambda x)=\lambda^\alpha H(x)$ for $\lambda\ge 0$.
Consider the fixed period problem
\be \left\{\matrix{\dot{x}(t)=JH^\prime(x(t)), \cr
     x(1)=x(0).         \cr }\right. \lb{2.22}\ee

Define
\be L_0^{\frac{\alpha}{\alpha-1}}(S^1,\R^{2n})
  =\{u\in L^{\frac{\alpha}{\alpha-1}}(S^1,\R^{2n})\,|\,\int_0^1udt=0\}. \lb{2.23}\ee
The corresponding Clarke-Ekeland dual action functional is defined by
\be \Phi(u)=\int_0^1\left(\frac{1}{2}Ju\cdot Mu+H^{\ast}(-Ju)\right)dt,
    \qquad \forall\;u\in L_0^{\frac{\alpha}{\alpha-1}}(S^1,\R^{2n}), \lb{2.24}\ee
where $Mu$ is defined by $\frac{d}{dt}Mu(t)=u(t)$ and $\int_0^1Mu(t)dt=0$,
$H^\ast$ is the Fenchel transform of $H$ defined above.

For any $\kappa\in\R$, we denote by
\be \Phi^{\kappa-}=\{u\in L_0^{\frac{\alpha}{\alpha-1}}(S^1,\R^{2n})\;|\;
             \Phi(u)<\kappa\}. \lb{2.25}\ee
Then as in P.218 of \cite{Eke3}, we define
\be c_i=\inf\{\delta\in\R\;|\: \hat I(\Phi^{\delta-})\ge i\},\lb{2.26}\ee
where $\hat I$ is the Fadell-Rabinowitz index given above. Then by Proposition 3
in P.218 of \cite{Eke3}, we have

{\bf Proposition 2.8.} {\it Every $c_i$ is a critical value of $\Phi$. If
$c_i=c_j$ for some $i<j$, then there are infinitely many geometrically
distinct closed characteristics on $\Sg$.}

As in Definition 2.1, we define the following

{\bf Definition 2.9.} {\it Suppose $u$ is a nonzero critical
point of $\Phi$, and $\Nn$ is an $S^1$-invariant
open neighborhood of $S^1\cdot u$ such that
$crit(\Phi)\cap(\Lambda(u)\cap \Nn)=S^1\cdot u$. Then
the $S^1$-critical modules of $S^1\cdot u$ is defined by
\bea C_{S^1,\; q}(\Phi, \;S^1\cdot u)
=H_{q}((\Lambda(u)\cap\Nn)_{S^1},\;
((\Lambda(u)\setminus S^1\cdot u)\cap\Nn)_{S^1}),\lb{2.27}
\eea
where $\Lambda(u)=\{w\in L_0^{\frac{\alpha}{\alpha-1}}(S^1,\R^{2n})\;|\;
\Phi(w)\le\Phi(u)\}$.}

Comparing with Theorem 4 in P.219 of \cite{Eke3}, we have the following

{\bf Proposition 2.10.} {\it For every $i\in\N$, there exists a point
$u\in L_0^{\frac{\alpha}{\alpha-1}}(S^1,\R^{2n})$ such that}
\bea
&& \Phi^\prime(u)=0,\quad \Phi(u)=c_i, \lb{2.28}\\
&& C_{S^1,\; 2(i-1)}(\Phi, \;S^1\cdot u)\neq 0. \lb{2.29}\eea

{\bf Proof.} By Lemma 8 in P.206 of \cite{Eke3}, we can use
Theorem 1.4.2 of \cite{Cha1} in the equivariant form  to obtain
\be H_{S^1,\,\ast}(\Phi^{c_i+\epsilon},\;\Phi^{c_i-\epsilon})
=\bigoplus_{\Phi(u)=c_i}C_{S^1,\; \ast}(\Phi, \;S^1\cdot u),\lb{2.30}\ee
for $\epsilon$  small enough such that the interval
$(c_i-\epsilon,\,c_i+\epsilon)$ contains no critical values of $\Phi$
except $c_i$.

Similar to P.431 of \cite{EkH1}, we have
\be H^{2(i-1)}((\Phi^{c_i+\epsilon})_{S^1},\,(\Phi^{c_i-\epsilon})_{S^1})
\mapright{q^\ast} H^{2(i-1)}((\Phi^{c_i+\epsilon})_{S^1} )
\mapright{p^\ast}H^{2(i-1)}((\Phi^{c_i-\epsilon})_{S^1}),
\lb{2.31}\ee
where $p$ and $q$ are natural inclusions. Denote by
$f: (\Phi^{c_i+\epsilon})_{S^1}\rightarrow CP^\infty$ a classifying map and let
$f^{\pm}=f|_{(\Phi^{c_i\pm\epsilon})_{S^1}}$. Then clearly each
$f^{\pm}: (\Phi^{c_i\pm\epsilon})_{S^1}\rightarrow CP^\infty$ is a classifying
map on $(\Phi^{c_i\pm\epsilon})_{S^1}$. Let $\eta \in H^2(CP^\infty)$ be the first
universal Chern class.

By definition of $c_i$, we have $\hat I(\Phi^{c_i-\epsilon})< i$, hence
$(f^-)^\ast(\eta^{i-1})=0$. Note that
$p^\ast(f^+)^\ast(\eta^{i-1})=(f^-)^\ast(\eta^{i-1})$.
Hence the exactness of (\ref{2.31}) yields a
$\sigma\in H^{2(i-1)}((\Phi^{c_i+\epsilon})_{S^1},\,(\Phi^{c_i-\epsilon})_{S^1})$
such that $q^\ast(\sigma)=(f^+)^\ast(\eta^{i-1})$.
Since $\hat I(\Phi^{c_i+\epsilon})\ge i$, we have $(f^+)^\ast(\eta^{i-1})\neq 0$.
Hence $\sigma\neq 0$, and then
$$H^{2(i-1)}_{S^1}(\Phi^{c_i+\epsilon},\Phi^{c_i-\epsilon})=
H^{2(i-1)}((\Phi^{c_i+\epsilon})_{S^1},\,(\Phi^{c_i-\epsilon})_{S^1})\neq 0. $$
Now the proposition follows from (\ref{2.30}) and the universal coefficient
theorem. \hfill\hb

{\bf Proposition 2.11.} {\it Suppose $u$ is the critical point of $\Phi$ found
in Proposition 2.10. Then we have
\be C_{S^1,\; 2(i-1)}(\Psi_a, \;S^1\cdot u_a)\cong C_{S^1,\; 2(i-1)}(\Phi,\;S^1\cdot u)\neq 0, \lb{2.32}\ee
where $\Psi_a$ is given by (\ref{2.3}) and $u_a\in L_0^2(S^1,\;\R^{2n})$
is its critical point corresponding to $u$ in the natural sense.}

{\bf Proof.} Fix this $u$, we  modify the function $H$ only in a small
neighborhood $\Omega$ of $0$ as in \cite{Eke1} so that the corresponding
orbit of $u$ does not enter $\Omega$ and the resulted function $\widetilde{H}$
satisfies similar properties as Definition 1 in P. 26
of \cite{Eke1} by just replacing $\frac{3}{2}$ there by $\alpha$.
Define the dual action functional
$\widetilde{\Phi}:L_0^{\frac{\alpha}{\alpha-1}}(S^1,\R^{2n})\to\R$ by
\be \widetilde{\Phi}(v)=\int_0^1\left(\frac{1}{2}Jv\cdot
   Mv+\widetilde{H}^{\ast}(-Jv)\right)dt. \lb{2.33}\ee
Clearly $\Phi$ and $\widetilde{\Phi}$ are $C^1$ close to each other,
thus by the continuity of critical modules (cf. Theorem 8.8 of \cite{MaW1} or
Theorem 1.5.6 in P.53 of \cite{Cha1}, which can be easily generalized to the
equivariant sense) for the $u$ in the proposition, we have
\be C_{S^1,\; \ast}(\Phi, \;S^1\cdot u)\cong C_{S^1,\; \ast}(\widetilde{\Phi},
    \;S^1\cdot u).\lb{2.34}\ee

Using a finite dimensional approximation as in Lemma 3.9 of \cite{Eke1},
we have
\be C_{S^1,\; \ast}(\widetilde{\Phi}, \;S^1\cdot u)
\cong H_\ast(\widetilde{\Lambda}(u)\cap D_\iota N(u),\;
    (\widetilde{\Lambda}(u)\setminus\{u\})\cap D_\iota N(u))^{\Z_m}, \lb{2.35}\ee
where $\widetilde{\Lambda}(u)=\{w\in L_0^{\frac{\alpha}{\alpha-1}}(S^1,\R^{2n})\;|\;
\widetilde{\Phi}(w)\le\widetilde{\Phi}(u)\}$ and $D_\iota N(u)$ is a
$\Z_m$-invariant finite dimensional disk transversal to $S^1\cdot u$ at $u$
(cf. Lemma 3.9 of \cite{WHL}), $m$ is the multiplicity of $u$.

By Lemma 3.9 of \cite{WHL}, we have
\be C_{S^1,\; \ast}(\Psi_a, \;S^1\cdot u_a)
\cong H_\ast(\Lambda_a(u_a)\cap D_\iota N(u_a),\;
    (\Lambda_a(u_a)\setminus\{u_a\})\cap D_\iota N(u_a))^{\Z_m}.\lb{2.36}\ee
By the construction of $H_a$ in \cite{WHL}, $H_a=\widetilde{H}$ in a
$L^\infty$-neighborhood of $S^1\cdot u$. We remark here that multiplying $H$ by a constant
will not affect the corresponding critical modules, i.e., the corresponding
critical orbits have isomorphic critical modules. Hence we can assume
$H_a=H$ in a $L^\infty$-neighborhood of $S^1\cdot u$ and then the above
conclusion holds. Hence $\Psi_a$ and $\widetilde{\Phi}$ coincide
in a $L^\infty$-neighborhood of $S^1\cdot u$. Note also by Lemma 3.9 of
\cite{Eke1}, the two finite dimensional approximations are actually the same.
Hence we have
\bea
&& H_\ast(\widetilde{\Lambda}(u)\cap D_\iota N(u),\;
   (\widetilde{\Lambda}(u)\setminus\{u\})\cap D_\iota N(u))^{\Z_m}\nn\\
&&\quad\cong H_\ast(\Lambda_a(u_a)\cap D_\iota N(u_a),\;
    (\Lambda_a(u_a)\setminus\{u_a\})\cap D_\iota N(u_a))^{\Z_m}.\lb{2.37}\eea
Now the proposition follows from Proposition 2.10 and (\ref{2.34})-(\ref{2.37}).
\hfill\hb

\setcounter{equation}{0}
\section{ A brief review on an index theory for symplectic paths}

In this section, we recall briefly an index theory for symplectic paths
developed by Y. Long and his coworkers.
All the details can be found in \cite{Lon4}.

As usual, the symplectic group $\Sp(2n)$ is defined by
$$ \Sp(2n) = \{M\in {\rm GL}(2n,\R)\,|\,M^TJM=J\}, $$
whose topology is induced from that of $\R^{4n^2}$. For $\tau>0$ we are interested
in paths in $\Sp(2n)$:
$$ \P_{\tau}(2n) = \{\ga\in C([0,\tau],\Sp(2n))\,|\,\ga(0)=I_{2n}\}, $$
which is equipped with the topology induced from that of $\Sp(2n)$. The
following real function was introduced in \cite{Lon3}:
$$ D_{\om}(M) = (-1)^{n-1}\ol{\om}^n\det(M-\om I_{2n}), \qquad
          \forall \om\in\U,\, M\in\Sp(2n). $$
Thus for any $\om\in\U$ the following codimension $1$ hypersurface in $\Sp(2n)$ is
defined in \cite{Lon3}:
$$ \Sp(2n)_{\om}^0 = \{M\in\Sp(2n)\,|\, D_{\om}(M)=0\}.  $$
For any $M\in \Sp(2n)_{\om}^0$, we define a co-orientation of $\Sp(2n)_{\om}^0$
at $M$ by the positive direction $\frac{d}{dt}Me^{t\ep J}|_{t=0}$ of
the path $Me^{t\ep J}$ with $0\le t\le 1$ and $\ep>0$ being sufficiently
small. Let
\bea
\Sp(2n)_{\om}^{\ast} &=& \Sp(2n)\bs \Sp(2n)_{\om}^0,   \nn\\
\P_{\tau,\om}^{\ast}(2n) &=&
      \{\ga\in\P_{\tau}(2n)\,|\,\ga(\tau)\in\Sp(2n)_{\om}^{\ast}\}, \nn\\
\P_{\tau,\om}^0(2n) &=& \P_{\tau}(2n)\bs  \P_{\tau,\om}^{\ast}(2n).  \nn\eea
For any two continuous arcs $\xi$ and $\eta:[0,\tau]\to\Sp(2n)$ with
$\xi(\tau)=\eta(0)$, it is defined as usual:
$$ \eta\ast\xi(t) = \left\{\matrix{
            \xi(2t), & \quad {\rm if}\;0\le t\le \tau/2, \cr
            \eta(2t-\tau), & \quad {\rm if}\; \tau/2\le t\le \tau. \cr}\right. $$
Given any two $2m_k\times 2m_k$ matrices of square block form
$M_k=\left(\matrix{A_k&B_k\cr
                                C_k&D_k\cr}\right)$ with $k=1, 2$,
as in \cite{Lon4}, the $\;\dm$-product of $M_1$ and $M_2$ is defined by
the following $2(m_1+m_2)\times 2(m_1+m_2)$ matrix $M_1\dm M_2$:
$$ M_1\dm M_2=\left(\matrix{A_1&  0&B_1&  0\cr
                               0&A_2&  0&B_2\cr
                             C_1&  0&D_1&  0\cr
                               0&C_2&  0&D_2\cr}\right). \nn$$  
Denote by $M^{\dm k}$ the $k$-fold $\dm$-product $M\dm\cdots\dm M$. Note
that the $\dm$-product of any two symplectic matrices is symplectic. For any two
paths $\ga_j\in\P_{\tau}(2n_j)$ with $j=0$ and $1$, let
$\ga_0\dm\ga_1(t)= \ga_0(t)\dm\ga_1(t)$ for all $t\in [0,\tau]$.

A special path $\xi_n$ is defined by
\be \xi_n(t) = \left(\matrix{2-\frac{t}{\tau} & 0 \cr
                                             0 &  (2-\frac{t}{\tau})^{-1}\cr}\right)^{\dm n}
        \qquad {\rm for}\;0\le t\le \tau.  \lb{3.1}\ee

{\bf Definition 3.1.} (cf. \cite{Lon3}, \cite{Lon4}) {\it For any $\om\in\U$ and
$M\in \Sp(2n)$, define
\be  \nu_{\om}(M)=\dim_{\C}\ker_{\C}(M - \om I_{2n}).  \lb{3.2}\ee
For any $\tau>0$ and $\ga\in \P_{\tau}(2n)$, define
\be  \nu_{\om}(\ga)= \nu_{\om}(\ga(\tau)).  \lb{3.3}\ee

If $\ga\in\P_{\tau,\om}^{\ast}(2n)$, define
\be i_{\om}(\ga) = [\Sp(2n)_{\om}^0: \ga\ast\xi_n],  \lb{3.4}\ee
where the right hand side of (\ref{3.4}) is the usual homotopy intersection
number, and the orientation of $\ga\ast\xi_n$ is its positive time direction under
homotopy with fixed end points.

If $\ga\in\P_{\tau,\om}^0(2n)$, we let $\mathcal{F}(\ga)$
be the set of all open neighborhoods of $\ga$ in $\P_{\tau}(2n)$, and define
\be i_{\om}(\ga) = \sup_{U\in\mathcal{F}(\ga)}\inf\{i_{\om}(\beta)\,|\,
                       \beta\in U\cap\P_{\tau,\om}^{\ast}(2n)\}.
               \lb{3.5}\ee
Then
$$ (i_{\om}(\ga), \nu_{\om}(\ga)) \in \Z\times \{0,1,\ldots,2n\}, $$
is called the index function of $\ga$ at $\om$. }

Note that when $\om=1$, this index theory was introduced by
C. Conley-E. Zehnder in \cite{CoZ1} for the non-degenerate case with $n\ge 2$,
Y. Long-E. Zehnder in \cite{LZe1} for the non-degenerate case with $n=1$,
and Y. Long in \cite{Lon1} and C. Viterbo in \cite{Vit2} independently for
the degenerate case. The case for general $\om\in\U$ was defined by Y. Long
in \cite{Lon3} in order to study the index iteration theory (cf. \cite{Lon4}
for more details and references).

For any symplectic path $\ga\in\P_{\tau}(2n)$ and $m\in\N$,  we
define its $m$-th iteration $\ga^m:[0,m\tau]\to\Sp(2n)$ by
\be \ga^m(t) = \ga(t-j\tau)\ga(\tau)^j, \qquad
  {\rm for}\quad j\tau\leq t\leq (j+1)\tau,\;j=0,1,\ldots,m-1.
     \lb{3.6}\ee
We still denote the extended path on $[0,+\infty)$ by $\ga$.

{\bf Definition 3.2.} (cf. \cite{Lon3}, \cite{Lon4}) {\it For any $\ga\in\P_{\tau}(2n)$,
we define
\be (i(\ga,m), \nu(\ga,m)) = (i_1(\ga^m), \nu_1(\ga^m)), \qquad \forall m\in\N.
   \lb{3.7}\ee
The mean index $\hat{i}(\ga,m)$ per $m\tau$ for $m\in\N$ is defined by
\be \hat{i}(\ga,m) = \lim_{k\to +\infty}\frac{i(\ga,mk)}{k}. \lb{3.8}\ee
For any $M\in\Sp(2n)$ and $\om\in\U$, the {\it splitting numbers} $S_M^{\pm}(\om)$
of $M$ at $\om$ are defined by
\be S_M^{\pm}(\om)
     = \lim_{\ep\to 0^+}i_{\om\exp(\pm\sqrt{-1}\ep)}(\ga) - i_{\om}(\ga),
   \lb{3.9}\ee
for any path $\ga\in\P_{\tau}(2n)$ satisfying $\ga(\tau)=M$.}

For a given path $\gamma\in {\cal P}_{\tau}(2n)$ we consider to deform
it to a new path $\eta$ in ${\cal P}_{\tau}(2n)$ so that
\begin{equation}
i_1(\gamma^m)=i_1(\eta^m),\quad \nu_1(\gamma^m)=\nu_1(\eta^m), \quad
         \forall m\in {\bf N}, \label{3.10}
\end{equation}
and that $(i_1(\eta^m),\nu_1(\eta^m))$ is easy enough to compute. This
leads to finding homotopies $\delta:[0,1]\times[0,\tau]\to {\rm Sp}(2n)$
starting from $\gamma$ in ${\cal P}_{\tau}(2n)$ and keeping the end
points of the homotopy always stay in a certain suitably chosen maximal
subset of ${\rm Sp}(2n)$ so that (\ref{3.10}) always holds. In fact,  this
set was first discovered in \cite{Lon3} as the path connected component
$\Omega^0(M)$ containing $M=\gamma(\tau)$ of the set
\begin{eqnarray}
  \Omega(M)=\{N\in{\rm Sp}(2n)\,&|&\,\sigma(N)\cap{\bf U}=\sigma(M)\cap{\bf U}\;
{\rm and}\;  \nonumber\\
 &&\qquad \nu_{\lambda}(N)=\nu_{\lambda}(M),\;\forall\,
\lambda\in\sigma(M)\cap{\bf U}\}. \label{3.11}
\end{eqnarray}
Here $\Omega^0(M)$ is called the {\it homotopy component} of $M$ in
${\rm Sp}(2n)$.

In \cite{Lon3} and \cite{Lon4}, the following symplectic matrices were introduced
as {\it basic normal forms}:
\begin{eqnarray}
D(\lambda)=\left(\matrix{\lm & 0\cr
         0  & \lm^{-1}\cr}\right), &\quad& \lm=\pm 2,\lb{3.12}\\
N_1(\lm,b) = \left(\matrix{\lm & b\cr
         0  & \lm\cr}\right), &\quad& \lm=\pm 1, b=\pm1, 0, \lb{3.13}\\
R(\th)=\left(\matrix{\cos\th & -\sin\th\cr
        \sin\th  & \cos\th\cr}\right), &\quad& \th\in (0,\pi)\cup(\pi,2\pi),
                     \lb{3.14}\\
N_2(\om,b)= \left(\matrix{R(\th) & b\cr
              0 & R(\th)\cr}\right), &\quad& \th\in (0,\pi)\cup(\pi,2\pi),
                     \lb{3.15}\end{eqnarray}
where $b=\left(\matrix{b_1 & b_2\cr
               b_3 & b_4\cr}\right)$ with  $b_i\in\R$ and  $b_2\not=b_3$.

Splitting numbers possess the following properties:

{\bf Lemma 3.3.} (cf. \cite{Lon3} and Lemma 9.1.5 of \cite{Lon4}) {\it Splitting
numbers $S_M^{\pm}(\om)$ are well defined, i.e., they are independent of the choice
of the path $\ga\in\P_\tau(2n)$ satisfying $\ga(\tau)=M$ appeared in (\ref{3.9}).
For $\om\in\U$ and $M\in\Sp(2n)$, splitting numbers $S_N^{\pm}(\om)$ are constant
for all $N\in\Om^0(M)$. }

{\bf Lemma 3.4.} (cf. \cite{Lon3}, Lemma 9.1.5 and List 9.1.12 of \cite{Lon4})
{\it For $M\in\Sp(2n)$ and $\om\in\U$, there hold
\begin{eqnarray}
S_M^{\pm}(\om) &=& 0, \qquad {\it if}\;\;\om\not\in\sg(M).  \lb{3.16}\\
S_{N_1(1,a)}^+(1) &=& \left\{\matrix{1, &\quad {\rm if}\;\; a\ge 0, \cr
0, &\quad {\rm if}\;\; a< 0. \cr}\right. \lb{3.17}\eea

For any $M_i\in\Sp(2n_i)$ with $i=0$ and $1$, there holds }
\be S^{\pm}_{M_0\dm M_1}(\om) = S^{\pm}_{M_0}(\om) + S^{\pm}_{M_1}(\om),
    \qquad \forall\;\om\in\U. \lb{3.18}\ee

We have the following

{\bf Theorem 3.5.} (cf. \cite{Lon3} and Theorem 1.8.10 of \cite{Lon4}) {\it For
any $M\in\Sp(2n)$, there is a path $f:[0,1]\to\Om^0(M)$ such that $f(0)=M$ and
\be f(1) = M_1\dm\cdots\dm M_s,  \lb{3.19}\ee
where each $M_i$ is a basic normal form listed in (\ref{3.12})-(\ref{3.15})
for $1\leq i\leq s$.}

Let $\Sigma\in\H(2n)$. Using notations in \S1,
for any closed characteristic $(\tau,y)$ and $m\in\N$, we define
its $m$-th iteration $y^m:\R/(m\tau\Z)\to\R^{2n}$ by
\be y^m(t) = y(t-j\tau), \qquad {\rm for}\quad j\tau\leq t\leq (j+1)\tau,
       \quad j=0,1,2,\ldots, m-1. \lb{3.20}\ee
Note that this coincide with that in \S2.
We still denote by $y$ its extension to $[0,+\infty)$.

We define via Definition 3.2 the following
\bea  S^+(y) &=& S_{\ga_y(\tau)}^+(1),  \lb{3.21}\\
  (i(y,m), \nu(y,m)) &=& (i(\ga_y,m), \nu(\ga_y,m)),  \lb{3.22}\\
   \hat{i}(y,m) &=& \hat{i}(\ga_y,m),  \lb{3.23}\eea
for all $m\in\N$, where $\ga_y$ is the associated symplectic path of $(\tau,y)$.
Then we have the following.

{\bf Theorem 3.6.} (cf. Lemma 1.1 of \cite{LoZ1}, Theorem 15.1.1 of \cite{Lon4}) {\it Suppose
$(\tau,y)$ is a closed characteristic on $\Sg$. Then we have
\be i(y^m)\equiv i(m\tau ,y)=i(y, m)-n,\quad \nu(y^m)\equiv\nu(m\tau, y)=\nu(y, m),
       \qquad \forall m\in\N, \lb{3.24}\ee
where $i(y^m)$ and $\nu(y^m)$ are the index and nullity
defined in \S2. }

The following is the precise iteration formula, which is due to Y. Long (cf. Theorem 8.3.1 and Corollary 8.3.2 of \cite{Lon4}).

{\bf Theorem 3.7.} {\it Let $\gamma\in\{\xi\in
C([0,\tau],Sp(2n))\mid \xi(0)=I\}$, Then there exists a path $f\in
C([0,1],\Omega^0(\gamma(\tau))$ such that $f(0)=\gamma(\tau)$ and
\bea f(1)=&&N_1(1,1)^{\diamond p_-} \diamond I_{2p_0}\diamond
N_1(1,-1)^{\diamond p_+}
\diamond N_1(-1,1)^{\diamond q_-} \diamond (-I_{2q_0})\diamond
N_1(-1,-1)^{\diamond q_+}\nn\\
&&\diamond R(\theta_1)\diamond\cdots\diamond R(\theta_r)
\diamond N_2(\omega_1, u_1)\diamond\cdots\diamond N_2(\omega_{r_*}, u_{r_*}) \nn\\
&&\diamond N_2(\lm_1, v_1)\diamond\cdots\diamond N_2(\lm_{r_0}, v_{r_0})
\diamond M_0 \lb{3.25}\eea
where $ N_2(\omega_j, u_j) $s are
non-trivial and   $ N_2(\lm_j, v_j)$s  are trivial basic normal
forms; $\sigma (M_0)\cap U=\emptyset$; $p_-$, $p_0$, $p_+$, $q_-$,
$q_0$, $q_+$, $r$, $r_*$ and $r_0$ are non-negative integers;
$\omega_j=e^{\sqrt{-1}\alpha_j}$, $
\lambda_j=e^{\sqrt{-1}\beta_j}$; $\theta_j$, $\alpha_j$, $\beta_j$
$\in (0, \pi)\cup (\pi, 2\pi)$; these integers and real numbers
are uniquely determined by $\gamma(\tau)$. Then using the
functions defined in (\ref{1.4}).
\bea i(\gamma, m)=&&m(i(\gamma,
1)+p_-+p_0-r)+2\sum_{j=1}^r E\left(\frac{m\theta_j}{2\pi}\right)-r
-p_--p_0\nn\\&&-\frac{1+(-1)^m}{2}(q_0+q_+)+2\left(
\sum_{j=1}^{r_*}\varphi\left(\frac{m\alpha_j}{2\pi}\right)-r_*\right).
\lb{3.26}\eea
\bea \nu(\gamma, m)=&&\nu(\gamma,
1)+\frac{1+(-1)^m}{2}(q_-+2q_0+q_+)+2(r+r_*+r_0)\nn\\
&&-2\left(\sum_{j=1}^{r}\varphi\left(\frac{m\theta_j}{2\pi}\right)+
\sum_{j=1}^{r_*}\varphi\left(\frac{m\alpha_j}{2\pi}\right)
+\sum_{j=1}^{r_0}\varphi\left(\frac{m\beta_j}{2\pi}\right)\right)\lb{3.27}\eea
\bea \hat i(\gamma, 1)=i(\gamma, 1)+p_-+p_0-r+\sum_{j=1}^r
\frac{\theta_j}{\pi}.\lb{3.28}\eea
Where $N_1(1, \pm 1)=
\left(\matrix{ 1 &\pm 1\cr 0 & 1\cr}\right)$, $N_1(-1, \pm 1)=
\left(\matrix{ -1 &\pm 1\cr 0 & -1\cr}\right)$,
$R(\theta)=\left(\matrix{\cos\th &
                  -\sin\th\cr\sin\th & \cos\th\cr}\right)$,
$ N_2(\omega, b)=\left(\matrix{R(\th) & b
                  \cr 0 & R(\th)\cr}\right)$    with some
$\th\in (0,\pi)\cup (\pi,2\pi)$ and $b=
\left(\matrix{ b_1 &b_2\cr b_3 & b_4\cr}\right)\in\R^{2\times2}$,
such that $(b_2-b_3)\sin\theta>0$, if $ N_2(\omega, b)$ is
trivial; $(b_2-b_3)\sin\theta<0$, if $ N_2(\omega, b)$ is
non-trivial. We have $i(\gamma, 1)$ is odd if $f(1)=N_1(1, 1)$, $I_2$,
$N_1(-1, 1)$, $-I_2$, $N_1(-1, -1)$ and $R(\theta)$; $i(\gamma, 1)$ is
even if $f(1)=N_1(1, -1)$ and $ N_2(\omega, b)$; $i(\gamma, 1)$ can be any
integer if $\sigma (f(1)) \cap \U=\emptyset$.}

\setcounter{equation}{0}
\section{Proof of the main theorems}

In this section, we give the proof of the main theorems  by using Morse theory
and the index iteration theory developed by Long and his coworkers.

As Definition 1.1 of \cite{LoZ1}, we define

{ \bf Definition  4.1.} For $\alpha\in(1,2)$, we define a map
$\varrho_n\colon\H(2n)\to\N\cup\{ +\infty\}$
\be \varrho_n(\Sg)
= \left\{\matrix{+\infty, & {\rm if\;\;}^\#\mathcal{V}(\Sigma,\alpha)=+\infty, \cr
\min\left\{[\frac{i(x,1) + 2S^+(x) - \nu(x,1)+n}{2}]\,
\left|\frac{}{}\right.\,(\tau,x)\in\mathcal{V}_\infty(\Sigma, \alpha)\right\},
 & {\rm if\;\;} ^\#\mathcal{V}(\Sigma, \alpha)<+\infty, \cr}\right.  \lb{4.1}\ee
where $\mathcal{V}(\Sigma,\alpha)$ and $\mathcal{V}_\infty(\Sigma,\alpha)$ are
variationally visible and infinite variationally visible sets respectively given
by Definition 1.4 of \cite{LoZ1} (cf. Definition 15.3.3 of \cite{Lon4}).

{\bf Proof of Theorem 1.1.}  By Theorem 1.1 of \cite{LoZ1} (cf. Theorem 15.4.3 of \cite{Lon4}), we have $^{\#}\T(\Sg)\ge \varrho_n(\Sg)\ge [\frac{n}{2}]+1$.  Hence  in order to prove Theorem 1.1, we only need to consider the case that $n$ being odd.
We prove Theorem 1.1 by contradiction, we assume $^{\#}\T(\Sg)=[\frac{n}{2}]+1$,
i.e., there are exactly  $[\frac{n}{2}]+1$  geometrically distinct
closed characteristics  $\{(\tau_j, y_j)\}_{1\le j\le [\frac{n}{2}]+1}$ on $\Sg$, where $n\in 2\N+1$.

Denote  by $\ga_j\equiv \gamma_{y_j}$ the associated
symplectic path of $(\tau_j,\,y_j)$ for $1\le j\le [\frac{n}{2}]+1$. Then by Lemma 1.3 of \cite{LoZ1}
(cf. Lemma 15.2.4 of \cite{Lon4}), there exist $P_j\in \Sp(2n)$ and $M_j\in \Sp(2n-2)$ such
that
\be \ga_j(\tau_j)=P_j^{-1}(N_1(1,\,1)\dm M_j)P_j, \quad 1\le j\le [\frac{n}{2}]+1.
   \lb{4.2}\ee
Using the common index jump theorem (Theorems 4.3 and 4.4 of
\cite{LoZ1}, Theorems 11.2.1 and 11.2.2 of \cite{Lon4}), we obtain some
$(T, m_1, m_2, \ldots,m_{[\frac{n}{2}]+1})\in\N^{[\frac{n}{2}]+2}$ such that the following hold by (11.2.6), (11.2.7) and
(11.2.26) of \cite{Lon4}:
\bea
i(y_j,\, 2m_j) &\ge& 2T-\frac{e(\gamma_j(\tau_j))}{2}, \lb{4.3}\\
i(y_j,\, 2m_j)+\nu(y_j,\, 2m_j) &\le& 2T+\frac{e(\gamma_j(\tau_j))}{2}-1, \lb{4.4}\\
i(y_j,\, 2m_j+1) &=& 2T+i(y_j,\,1). \lb{4.5}\\
i(y_j,\, 2m_j-1)+\nu(y_j,\, 2m_j-1)
 &=& 2T-(i(y_j,\,1)+2S^+_{\gamma_j(\tau_j)}(1)-\nu(y_j, 1)),\lb{4.6}
\eea
where $e(M)$ is the total algebraic multiplicity of all eigenvalues of $M$ on $\U$. 

By Corollary 1.2 of \cite{LoZ1} (cf. Corollary 15.1.4 of \cite{Lon4}), we
have $i(y_j,\,1)\ge n$ for $1\le j\le [\frac{n}{2}]+1$. Note that
$e(\gamma_j(\tau_j))\le 2n$ for $1\le j\le [\frac{n}{2}]+1$. Hence Theorem 2.3 of \cite{LoZ1}  (cf. Theorem 10.2.4 of \cite{Lon4}) yields
\bea i(y_j,\, m)+\nu(y_j,\, m)
&\le&  i(y_j, m+1)-i(y_j, 1)+\frac{e(\gamma_j(\tau_j))}{2}-1\nn\\
&\le& i(y_j, m+1)-1. \quad \forall m\in\N,\;1\le j\le [\frac{n}{2}]+1.\lb{4.7}\eea

By Theorem 3.7, the matrix $M_j$ can be connected within $\Omega^0(M_j)$ to $N_1(1,1)^{\diamond p_{j_-}} \diamond I_{2p_{j_0}}\diamond
N_1(1,-1)^{\diamond p_{j_+}}
\diamond M_j^\prime$, where $p_{j_-}, \,p_{j_0},\,p_{j_+}\in\N_0$ and $1\notin\sigma(M_j^\prime)$ for $1\le j\le [\frac{n}{2}]+1$.
By Lemma 3.4, we have
\bea
&& 2S^+_{\gamma_j(\tau_j)}(1)-\nu(y_j,\,1) \nn\\
&&\qquad = 2S^+_{N_1(1,\,1)}(1)-\nu_1(N_1(1,\,1))
   +2S^+_{M_j}(1)-\nu_1(M_j)\nn\\
&&\qquad = 1 + p_{j_-}-p_{j_+}, \qquad 1\le j\le [\frac{n}{2}]+1.\lb{4.8}\eea
Note that by (\ref{4.1}) and (\ref{4.8}), we have $\varrho_n(\Sg)\ge [\frac{n+1}{2}]+1$
if there is no closed characteristic $(\tau_j,\,y_j)$  on $\Sg$ satisfies $p_{j_+}=n-1$ together with $i(y_j,\, 1)=n$. Thus in oder to prove Theorem 1.1, it is sufficient to consider the case that there is some $(\tau_j,\,y_j)$ satisfies  $p_{j_+}=n-1$ together with $i(y_j,\, 1)=n$ for some $j\in\{1, 2,\ldots,[\frac{n}{2}]+1\}$, i.e.,  the matrix $M_j$ can be connected within $\Omega^0(M_j)$ to $N_1(1,-1)^{\diamond (n-1)}$.  Without loss of generality, we may assume $i(y_j,\,1)=n$ and  $M_j$ can be connected within $\Omega^0(M_j)$ to $N_1(1,-1)^{\diamond (n-1)}$ whenever $1\le j\le K$,
 where $K$ is some integer in $[1.\, [\frac{n}{2}]+1]$.
 
By Theorem 3.6, (\ref{4.7}) and (\ref{4.8}), (\ref{4.3})-(\ref{4.6}) become
\bea
i(y_j^{2m_j}) &\ge& 2T-2n, \quad 1\le j\le [\frac{n}{2}]+1,\lb{4.9}\\
i(y_j^{2m_j})+\nu(y_j^{2m_j})-1 &\le& 2T-2, \quad 1\le j\le [\frac{n}{2}]+1,  \lb{4.10}\\
i(y_j^{2m_j+m}) &\ge& 2T, \quad\forall\; m\ge 1, \quad 1\le j\le [\frac{n}{2}]+1, \lb{4.11}\\
i(y_j^{2m_j-1})+\nu(y_j^{2m_j-1})-1 &=& 2T-n-3,\quad 1\le j\le K, \lb{4.12}\\
i(y_j^{2m_j-m})+\nu(y_j^{2m_j-m})-1 &<& 2T-n-3,\quad\forall\; m\ge 2,\; 1\le j\le K,\lb{4.13}\\
i(y_j^{2m_j-m})+\nu(y_j^{2m_j-m})-1 &<& 2T-n-3,\quad\forall\; m\ge 1, \;K< j\le[\frac{n}{2}]+1.\lb{4.14}
\eea
Thus by Propositions 2.10 and 2.11,  we can find $(j_k,\,l_{j_k})_{1\le k<\infty}$ such that
\bea &&\Phi^\prime(u_{j_k}^{l_{j_k}})=0,\quad \Phi(u_{j_k}^{l_{j_k}})=c_{T+1-k},
\qquad C_{S^1,\; 2T-2k}(\Psi_a, \;S^1\cdot u_{j_k}^{l_{j_k}})\neq 0,
\lb{4.15}
\eea
where we denote also by $u_{j_k}^{l_{j_k}}$ the corresponding
critical points of $\Phi$ and which will not be confused.
Tlhus by Propositions 2.3, 2.8  and (\ref{4.9})-(\ref{4.15}), we have 
$(j_k,\,l_{j_k})=(j_k,\,2m_{j_k})$ for $1\le k\le [\frac{n}{2}]+1$ and these $j_k$s
are pairwise distinct.

{\bf Claim. } {\it We have $1\le j_1\le K$.}

We prove the claim by contradiction.  Suppose $j_1\notin [1,\,K]$.
Then by  Theorem 3.7, Proposition 2.3,  2.5 and (\ref{4.10})-(\ref{4.15}), we have  
$C_{S^1,\; 2T-n-3}(\Psi_a, \;S^1\cdot u_j^{2m_j-1})=0$ for $1\le j\le [\frac{n}{2}]+1$.
In fact, we have  $C_{S^1,\; 2T-n-3}(\Psi_a, \;S^1\cdot u_j^{2m_j-1})=0$ for $K< j\le [\frac{n}{2}]+1$  by (\ref{4.14}).  Since $j_1\notin [1,\,K]$, we have  $C_{S^1,\; 2T-2k}(\Psi_a, \;S^1\cdot u_{j_k}^{2m_{j_k}})\neq 0$ for  $j_k\in[1,\,K]$ and $k\neq 1$ by (\ref{4.15}). Thus it follows from (\ref{4.10}), Propositions 2.3 and 2.6 that   $C_{S^1,\; 2T-2}(\Psi_a, \;S^1\cdot u_{j}^{2m_{j}})= 0$ for $1\le j\le K$, and then  $k_{\nu(y_j)-1}(y_j)=0$  for $1\le j\le K$ by Theorem 3.7 and Proposition 2.5. Thus by Proposition 2.5 and (\ref{4.12}), we have
$C_{S^1,\; 2T-n-3}(\Psi_a, \;S^1\cdot u_j^{2m_j-1})=0$ for $1\le j\le K$.
This together with  (\ref{4.13})-(\ref{4.15}) implies  $M_{2T-n-3}=0$. This contradict to
Theorem 2.7. Hence the claim holds.
 
 By the proof of Theorem 1.4 of \cite{LoZ1} (cf. Theorem 15.5.3 of \cite{Lon4}), we can find  another tuple $(T^\prime, m_1^\prime, m^\prime_2, \ldots,m^\prime_{[\frac{n}{2}]+1})\in\N^{[\frac{n}{2}]+2}$ such that the following hold \bea
i(y_j^{2m^\prime_j}) &\ge& 2T^\prime-2n, \quad 1\le j\le [\frac{n}{2}]+1,\lb{4.16}\\
i(y_j^{2m^\prime_j})+\nu(y_j^{2m^\prime_j})-1 &\le& 2T^\prime-2, \quad 1\le j\le [\frac{n}{2}]+1,  \lb{4.17}\\
i(y_j^{2m^\prime_j+m}) &\ge& 2T^\prime, \quad\forall\; m\ge 1, \quad 1\le j\le [\frac{n}{2}]+1, \lb{4.18}\\
i(y_j^{2m^\prime_j-1})+\nu(y_j^{2m^\prime_j-1})-1 &=& 2T^\prime-n-3,\quad 1\le j\le K, \lb{4.19}\\
i(y_j^{2m^\prime_j-m})+\nu(y_j^{2m^\prime_j-m})-1 &<& 2T^\prime-n-3,\quad\forall\; m\ge 2,\; 1\le j\le K,\lb{4.20}\\
i(y_j^{2m^\prime_j-m})+\nu(y_j^{2m^\prime_j-m})-1 &<& 2T^\prime-n-3,\quad\forall\; m\ge 1, \;K< j\le[\frac{n}{2}]+1.\lb{4.21}
\eea  
As above we can find $(i_k,\,l_{i_k})_{1\le k<\infty}$ such that
\bea &&\Phi^\prime(u_{i_k}^{l_{i_k}})=0,\quad \Phi(u_{i_k}^{l_{i_k}})=c_{T^\prime+1-k},
\qquad C_{S^1,\; 2T^\prime-2k}(\Psi_a, \;S^1\cdot u_{i_k}^{l_{i_k}})\neq 0,
\lb{4.22}
\eea
together with
$(i_k,\,l_{i_k})=(i_k,\,2m^\prime_{i_k})$ for $1\le k\le [\frac{n}{2}]+1$ and these $i_k$s
are pairwise distinct. 

By the proof of Theorem 1.4 of \cite{LoZ1} (cf. Theorem 15.5.3 of \cite{Lon4}), we can  require  the tuple  $(T^\prime, m_1^\prime, m_2^\prime,\ldots, m_{[\frac{n}{2}]+1}^\prime)$ further satisfies $i_1\neq j_1$. In fact, otherwise we always have the same elliptic closed lcharacteristic $(\tau_{j_1}, y_{j_1})$ for all $a\in A(v)$, where we use notations as in Theorem 1.4 of \cite{LoZ1} (cf. Theorem 15.5.3 of \cite{Lon4}). Then the  proof of Theorem 1.4 of \cite{LoZ1} (cf. Theorem 15.5.3 of \cite{Lon4}) yields a contradiction. 

On the other hand, by Theorem 3.7, Propositions 2.3, 2.5 and the above claim,
we have $C_{S^1,\; 2T^\prime-2}(\Psi_a, \;S^1\cdot u_{j_1}^{2m_{j_1}^\prime})\neq0$.
Thus it follows from Propositions 2.3, 2.6 and (\ref{4.22}) that $j_1\notin\{i_2, i_3,\ldots,i_{[\frac{n}{2}]+1}\}$.
This contradict to the assumption
 that there are $[\frac{n}{2}]+1$ geometrically distinct closed characteristics on $\Sg$.
 The proof of Theorem 1.1 is complete.\hfill\hb

{\bf Proof of Theorem 1.2.}  Denote the closed characteristics  on $\Sg$ by $\{(\tau_j, y_j)\}_{1\le j\le q}$.
Denote  by $\ga_j\equiv \gamma_{y_j}$ the associated
symplectic path of $(\tau_j,\,y_j)$ for $1\le j\le q$. Then by Lemma 1.3 of \cite{LoZ1} (cf. Lemma
15.2.4 of \cite{Lon4}), there exist $P_j\in \Sp(2n)$ and $M_j\in \Sp(2n-2)$ such
that
\be \ga_j(\tau_j)=P_j^{-1}(N_1(1,\,1)\dm M_j)P_j, \quad 1\le j\le q.
   \lb{4.23}\ee
Using the common index jump theorem (Theorems 4.3 and 4.4 of
\cite{LoZ1}, Theorems 11.2.1 and 11.2.2 of \cite{Lon4}), we obtain some
$(T, m_1, m_2, \ldots,m_q)\in\N^{q+1}$ such that the following hold by (11.2.6), (11.2.7) and
(11.2.26) of \cite{Lon4}:
\bea
i(y_j,\, 2m_j) &\ge& 2T-\frac{e(\gamma_j(\tau_j))}{2}, \lb{4.24}\\
i(y_j,\, 2m_j)+\nu(y_j,\, 2m_j) &\le& 2T+\frac{e(\gamma_j(\tau_j))}{2}-1, \lb{4.25}\\
i(y_j,\, 2m_j+1) &=& 2T+i(y_j,\,1). \lb{4.26}\\
i(y_j,\, 2m_j-1)+\nu(y_j,\, 2m_j-1)
 &=& 2T-(i(y_j,\,1)+2S^+_{\gamma_j(\tau_j)}(1)-\nu(y_j, 1)),\lb{4.27}
\eea
for $1\le j\le q$.
As in the proof of Theorem 1.1, we can find $(j_k,\,l_{j_k})_{1\le k<\infty}$ such that
\bea &&\Phi^\prime(u_{j_k}^{l_{j_k}})=0,\quad \Phi(u_{j_k}^{l_{j_k}})=c_{T+1-k},
\qquad C_{S^1,\; 2T-2k}(\Psi_a, \;S^1\cdot u_{j_k}^{l_{j_k}})\neq 0,
\lb{4.28}
\eea
and
$(j_k,\,l_{j_k})=(j_k,\,2m_{j_k})$ for $1\le k\le [\frac{n}{2}]+1$ and these $j_k$s
are pairwise distinct.

Note that by Theorem 3.6, (\ref{4.24}) and (\ref{4.25}),  a hyperbolic closed characteristic  $(\tau_j,\,y_j)$ must satisfy  $i(y_j^{2m_j})=2T-n-1$. Thus by Proposition 2.3 and (\ref{4.28}),
the closed characteristics $(\tau_{j_k},\,y_{j_k})$ for $1\le k\le [\frac{n}{2}]+1$ must be non-hyperbolic if $n$ is even.

It remains to consider the case that $n$ being odd. As above 
the closed characteristics $(\tau_{j_k},\,y_{j_k})$ for $1\le k\le [\frac{n}{2}]$  must be non-hyperbolic. We have the following two cases.

{\bf Case 1.} {\it  We have $j_1=j_{[\frac{n}{2}]+2}$.}

In this case, as in the proof of Theorem 1.1,  we must have $i(y_{j_1},\,1)=n$ and  $M_{j_1}$ can be connected within $\Omega^0(M_{j_1})$ to $N_1(1,-1)^{\diamond (n-1)}$ .
Hence by Theorem 3.7, we have $\hat i(y_{j_1})\in\Q$.  Then by the proof of Theorem 1.3 of \cite{LoZ1} (cf, Theorem 15.5.2 of \cite{Lon4}), we have $\hat i(y_{j_{[\frac{n}{2}]+1}})\notin\Q$. Thus $(\tau_{j_{[\frac{n}{2}]+1}},\,y_{j_{[\frac{n}{2}]+1}})$ is non-hyperbolic. 

 {\bf Case 2.} {\it  We have $j_1\neq j_{[\frac{n}{2}]+2}$.}
 
 By Propositions 2.3-2.6,  (\ref{4.24})-(\ref{4.27}), (\ref{4.7}) and (\ref{4.8}), we have $j_{[\frac{n}{2}]+2}\neq j_k$ for $2\le k\le [\frac{n}{2}]+1$. Hence the closed characteristics $(\tau_{j_k},\,y_{j_k})$ for $1\le k\le [\frac{n}{2}]$ and $k=[\frac{n}{2}]+2$ must be geometrically distinct non-hyperbolic closed characteristics. 
 
The proof of Theorem 1.2 is complete.\hfill\hb

\noindent {\bf Acknowledgements.} I would like to sincerely thank my
Ph. D. thesis advisor, Professor Yiming Long, for introducing me to Hamiltonian
dynamics and for his valuable help and encouragement during my research. I would like to say how enjoyable it is to work with him.

\bibliographystyle{abbrv}

\medskip

\end{document}